\documentclass[12pt]{amsart}
\usepackage[dviwindo]{graphicx}
\usepackage{amsmath,amsthm,latexsym,amssymb,amscd,verbatim,pifont,accents,yfonts,anysize}
\usepackage{pstricks,pst-eps}
\usepackage{url}
\usepackage{stmaryrd}

\usepackage{bm}
\usepackage[T1]{fontenc}

\usepackage{libertine}

 \theoremstyle{plain}

\input diagxy
   \theoremstyle{definition}

\newcommand{\ga}{\alpha}

\newcommand{\gk}{\kappa}
\newcommand{\gf}{\varphi}
\newcommand{\gS}{\Sigma}

\newcommand{\go}{\omega}

\newcommand{\y}{\wedge}

\newcommand{\e}{\exists\,}
\newcommand{\f}{\forall\,}

\newcommand{\set}[1]{ \{ #1 \}}

\newcommand{\vx}{\vec{x}}

\newcommand{\vb}{\vec{b}}

\renewcommand{\qed}{\vspace{0.3cm}\hfill\ding{112}}

\newcommand{\rl}{\blacktriangleright}

\newcommand{\il}{\blacktriangleleft}

\newcommand{\Ci}[1]{\mathring{#1}}

\newcommand{\C}{\mathfrak{C}}

\newcommand{\A}{\mathfrak{A}}

\numberwithin{Thm}{section} \numberwithin{equation}{section} 
\begin{document}

\title[About the review]{About the review in Mathematical Reviews of my paper:\\
The two-cardinal problem for languages of arbitrary cardinality\\
The Journal of Symbolic Logic {\bf 75}, Number 3, Sept., 2010, pp.  785-801\\
doi:10.2178/jsl/1278682200}
\author[Luis Villegas]{Luis Miguel Villegas Silva}

 \address{Departamento de Matem\'aticas, Universidad Aut\'onoma
Me\-tro\-po\-li\-ta\-na Iztapalapa, Av. San Rafael Atlixco 186, Col. Vicentina,
Iztapalapa, 09340 D.F., M\'exico\\ {\bf e-mail: lmvs@\-xanum.\-uam.\-mx}}

\maketitle

In this note I make some comments on the review
MR2723767 (2011m:03064)
appeared in Mathematical Reviews. I can not provide access here to the review or even to my paper, it is copyright material.

\section{The review}
\begin{enumerate}
\item  Paragraph five of the review:  "The paper observes that, in $L$,there is a coarse...."\\
In the paper (page 788) I give the axioms of a $(\mu,1)$-coarse morass according to Jensen  in terms of pairs of primitive recursive closed ordinals, and I mention that " A proof of the existence of a $(\mu,1)$-Coarse morass in $L$ can be extracted from [Dev84] for pairs of adequate ordinals"...

An adequate ordinal is  admissible or the limit of  admissible ordinals (Devlin p. 339 at the Bottom). In particular, they are primitive recursive closed.  The proof of the existence of a full $(\mu,1)$-morass starts in page 344 in Devlin's book. Actually what Devlin really construct is a $(\omega_1,1)$-morass, but the proof works perfectly good for any regular cardinal other that $\go_1$.  It is just a matter of following this proof, to corroborate that he constructs, in particular, a $(\go_1,1)$-coarse morass.   In item (b) I should write closed on $\sup(S_\ga)$, not in $\mu^+$ as the reviewer claims.

\item  Again Paragraph five of the review: "For example, on page  793, lines 1-3, a function is defined and stated to be $\gS_1(\set{\ga_1})$ but it is not....."

The constructions begins in page 792, not in page 793 as the reviewer wrote.  I claim that $S_{\ga_\nu}\cap\nu$  is $\gS_1(\set{\ga_\nu)}$ for every $\nu\in S^1$, which follows from any construction of a $(\mu,1)$-morass, because this is necessary to succeed on building such a morass.   For other  proofs of the existence of a $(\mu,1)$-morass see:  L. Stanley, {\em A short course on gap-one morasses with a review of the fine structure of $L$},  in: Surveys on set theory A. Mathias (Ed.), London Math. Soc. Lecture Notes Series \# 87, 1983,  pp. 197-244, or  P. Welch, {\em $\gS^*$-fine structure}, in A. Kanamori, M. Foreman (Eds), Handbook of Set Theory, Springer-Verlag, pp. 657-736, or the reference [Don81] in the paper.

What is important to us: I am not  claiming that the morass maps are $\gS_1$-preserving for a language which expands  LST, I only need the above mentioned fact that $S_{\ga_\nu}\cap\nu$  is $\gS_1(\set{\ga_\nu)}$ for every $\nu\in S^1$.  I also use that $<_L$ is $\gS_1$-definable to build the functions $h_{\ga_\nu}$. The sets $B_\nu$ appear in $P(L_{\ga_\nu})$. Indeed what we want is to enumerate those sets in a $\gS_1(\set{\ga_\nu})$-fashion. Once we have this enumeration, we can appeal to the morass maps and get the desired preservation.

\item Paragraph 6. Indeed Lemma 5.4 as stated is wrong,  we have to require that the $\vx$ belong to $U$.  But we can take  this lemma off the paper, it is not necessary in what follows.

    The proof of Lemma 5.6 is unnecessarily  complicated, and we do not need Lema 5.4. Let me provide a clearer proof.

    {\bf Proof of Lemma 5.6.} We keep  the given proof until the bottom of page 791. We have to show that $c$ is transcendent in $\C$ over $B\cup U^{\C}$.  First we prove that $c$ is transcendent in $\C$ over $B$. Otherwise there exists a formula $\gf(v,\vb)$, $\vb\in B$ in the complete type of $c$ in $\C$ over $B$ such that
    \begin{align*}
    (\C,c,\Ci{B})\models\gf(c,\vb)\\
    \intertext{and}
    (\C,\Ci{B})\models\neg\coprod v\gf(v,\vb).
    \end{align*}
    From the last assertion we get
    $$
    (\C,\Ci{B})\models\e u\f v\rl u\neg\gf(v,\vb),
    $$
    hence $\neg\gf(v,\vb)$ would appear in $\gS_c(v)$ by construction of this set (page 791), thus $(\C,c,\Ci{B})\models\neg\gf(c,\vb)$, a contradiction.

    Now we show that $c$ is transcendent in $\C$ over $U^{\C}$. If this is not the case,  as above, we find a formula $\gf(v,\vx)$ with $\vx\in U^{\C}$, such that
    \begin{align*}
    (\C,c,\vx)_{\vx\in U^{\C}}\models\gf(c,\vx)\\
    \intertext{and}
    (\C,\vx)_{\vx\in U^{\C}}\models\neg\coprod v\gf(v,\vx)
    \end{align*}

    Then
    $$
    (\C,c,\vx)_{\vx\in U^{\C}}\models\e\vx(U(\vx)\y\gf(c,\vx)\y\neg\coprod v\gf(v,\vx))
    $$

    The formula at the right has only $c$ as parameter, so it belongs to the complete type of $c$ in $\C$ over the empty set (or over $B$). Therefore
    $$
    \C\models\coprod w\e u\e\vx(U(\vx)\y\gf(w,\vx)\y\neg\coprod v\gf(v,\vx))
    $$
    The elements $\vx$ belong to $U$, which is lineraly ordered without maximum, so we can find $u\in U$ with
    $$
    \C\models\coprod w\e u\e\vx\il u(U(\vx)\y\gf(w,\vx)\y\neg\coprod v\gf(v,\vx))
    $$
    This is a contradiction:  the formula at the right has no parameters at all, and if $\A$ is the original structure of type $(\gk^+,\gk)$, we get
    $$
    \A\models\coprod w\e u\e\vx\il u(U(\vx)\y\gf(w,\vx)\y\neg\coprod v\gf(v,\vx))
    $$
    because of $\A\equiv\C$.  Then
    \begin{equation*}
    \A\models \f r\e w\rl r\e u\e\vx\il u(U(\vx)\y\gf(w,r)\y\neg\coprod v\gf(v,\vx))
    \end{equation*}
    For each $r\in A$ we find $u,\vx\in U$. We have available $\gk^+$ such $r$'s and only $\gk$ $u,\vx$'s, then there exist $u,\vx\in U$ such that
    $$
    \C\models\e u\e\vx\il u\coprod w(U(\vx)\y\gf(w,x)\y\neg\coprod v\gf(v,\vx))
    $$
    hence
    $$
    \C\models\e u\e\vx\il u(U(\vx)\y\coprod w\gf(w,\vx)\y\neg\coprod w\gf(w,\vx))
    $$
    which is clearly contradictory. \qed

    \item Last paragraph in the review:  "Let me add that the argument..."; it is often necessary to cite results from other researcher  or own results, and it is not always possible to anticipate any future developments.
        \end{enumerate}

\section{Personal remarks}

\end{document}